\documentclass[
11pt,%
tightenlines,%
twoside,%
onecolumn,%
nofloats,%
nobibnotes,%
nofootinbib,%
superscriptaddress,%
noshowpacs,%
centertags]%
{revtex4}
\usepackage{ljmarxiv}
\usepackage{dsfont}
\usepackage{placeins}

\newtheorem{cor}{Corollary}
\newtheorem{lem}{Lemma}

\newtheorem{que}{Question}
\newtheorem{rem}{Remark}
\newtheorem{prob}{Problem}
\newtheorem{conj}{Conjecture}
\theoremstyle{definition}
\newtheorem{example}[equation]{Example}

\newcounter {own}
\def\theown {\thesection       .\arabic{own}}

\newenvironment{pf}[1][]{%
 \vskip 3mm
 \noindent
 \ifthenelse{\equal{#1}{}}%
  {{\slshape Proof. }}%
  {{\slshape #1.} }%
 }%
{\qed\bigskip}

\newcounter{alphabet}

\newenvironment{Thm}[1][]{\refstepcounter{alphabet}%
\bigskip%
\noindent%
{\bf Theorem \Alph{alphabet}}%
\ifthenelse{\equal{#1}{}}{}{ (#1)}%
{\bf .} \itshape}{\vskip 8pt}

\newcommand{\IN}{{\mathbb N}}
\newcommand{\IC}{{\mathbb C}}
\newcommand{\ID}{{\mathbb D}}

\def\be{\begin{equation}}
\def\ee{\end{equation}}

\newcommand{\bee}{\begin{enumerate}}
\newcommand{\eee}{\end{enumerate}}

\newcommand{\blem}{\begin{lem}}
\newcommand{\elem}{\end{lem}}
\newcommand{\bthm}{\begin{theorem}}
\newcommand{\ethm}{\end{theorem}}
\newcommand{\bcor}{\begin{cor}}
\newcommand{\ecor}{\end{cor}}
\newcommand{\beg}{\begin{example}}
\newcommand{\eeg}{\end{example}}
\newcommand{\bques}{\begin{que}}
\newcommand{\eques}{\end{que}}
\newcommand{\begs}{\begin{examples}}
\newcommand{\eegs}{\end{examples}}
\newcommand{\bdefe}{\begin{defin}}
\newcommand{\edefe}{\end{defin}}
\newcommand{\bprob}{\begin{prob}}
\newcommand{\eprob}{\end{prob}}
\newcommand{\bei}{\begin{itemize}}
\newcommand{\eei}{\end{itemize}}

\newcommand{\bcon}{\begin{conj}}
\newcommand{\econ}{\end{conj}}
\newcommand{\bcons}{\begin{conjs}}
\newcommand{\econs}{\end{conjs}}
\newcommand{\bprop}{\begin{propo}}
\newcommand{\eprop}{\end{propo}}
\newcommand{\br}{\begin{rem}}
\newcommand{\er}{\end{rem}}
\newcommand{\brs}{\begin{rems}}
\newcommand{\ers}{\end{rems}}
\newcommand{\bo}{\begin{obser}}
\newcommand{\eo}{\end{obser}}
\newcommand{\bos}{\begin{obsers}}
\newcommand{\eos}{\end{obsers}}
\newcommand{\bpf}{\begin{pf}}
\newcommand{\epf}{\end{pf}}
\newcommand{\ba}{\begin{array}}
\newcommand{\ea}{\end{array}}
\newcommand{\beq}{\begin{eqnarray}}
\newcommand{\beqq}{\begin{eqnarray*}}
\newcommand{\eeq}{\end{eqnarray}}
\newcommand{\eeqq}{\end{eqnarray*}}

\newcommand{\ra}{\rightarrow}

\newcommand{\ds}{\displaystyle}

\begin{document}

\titlerunning{Bohr--Rogosinski inequalities} 
\authorrunning{Seraj A.~Alkhaleefah et al.} 

\title{Bohr--Rogosinski inequalities
for bounded analytic functions}

\author{\firstname{Seraj A.}~\surname{Alkhaleefah}}
\email[E-mail: ]{s.alkhaleefah@gmail.com}
\affiliation{N.I. Lobachevskii Institute of Mathematics and Mechanics, Kazan (Volga Region) Federal University, Kremlevskaya ul. 18, Kazan, Tatarstan, 420008 Russia}

\author{\firstname{Ilgiz R.}~\surname{Kayumov}}
\email[E-mail: ]{ikayumov@gmail.com}
\affiliation{N.I. Lobachevskii Institute of Mathematics and Mechanics, Kazan (Volga Region) Federal University, Kremlevskaya ul. 18, Kazan, Tatarstan, 420008 Russia}

\author{\firstname{Saminathan}~\surname{Ponnusamy}}
\email[E-mail: ]{samy@iitm.ac.in}
\affiliation{Department of Mathematics,
Indian Institute of Technology Madras,
Chennai-600 036, India}
%




\begin{abstract}
In this paper we first consider another version of the Rogosinski inequality for analytic functions
$f(z)=\sum_{n=0}^\infty a_nz^n$ in the unit disk $|z| < 1$, in which we replace the coefficients $a_n$ $(n= 0,1,\ldots ,N)$ of
the power series by the derivatives $f^{(n)}(z)/n!$ $(n= 0,1,\ldots ,N)$. Secondly, we
obtain improved versions of the classical Bohr inequality and Bohr's inequality for the harmonic mappings
of the form $f = h + \overline{g}$, where the analytic part $h$ is bounded by $1$ and
that $|g'(z)| \le k|h'(z)|$ in $|z| < 1$ and for some $k \in [0,1]$.
\end{abstract}

\maketitle


\section{Introduction}
Let $\ID=\{z\in \IC:\, |z| < 1 \}$  denote the open unit disk, and
${\mathcal A}$ denote the space of analytic functions in $\ID$ with
the topology of uniform convergence on compact sets. Define
${\mathcal B}=\{f\in {\mathcal A}:\, \mbox{$|f(z)|<1$ in $\ID$}\}$.
Then the Bohr radius is the largest number $r>0$ such that if $f\in
{\mathcal B} $ has the power series expansion
$f(z)=\sum_{n=0}^\infty a_n z^n$, then
$\sum_{n=0}^{\infty}|a_n|\,|z|^n \leq 1 ~\mbox{ for }~ |z|\leq r
$
which is called the classical Bohr inequality for the family
${\mathcal B}$. Rogosinski radius is the largest number $r>0$ such
that, under the previous assumptions, $\left|S_{N}(z) \right|<1
~\mbox{ for }~ |z|<r,
$
where $S_{N}(z)=\sum_{n=0}^{N}a_nz^n$ $(N\geq 0)$ denote the partial sums of $f$. This inequality is called
the classical Rogosinski inequality for the family ${\mathcal B}$.

 If $\textbf{B}$ and $\textbf{R}$ denote the Bohr radius and the
Rogosinski radius, respectively, then because $\left |S_{N}(z)\right |\leq \sum_{n=0}^{N}|a_n|\,|z|^n\leq \sum_{n=0}^{\infty}|a_n|\,|z|^n,$
it is clear that $\textbf{B}\leq \textbf{R}$. In fact the following two classical results are well-known.

\begin{Thm} Suppose that $f\in {\mathcal B} $. Then we have
{\rm $\textbf{B}=1/3$, and  {\rm (see Rogosinski \cite{Rogo-9} and also \cite{LanDai-86,SchSz-25})} $\textbf{R}=1/2$.}
\end{Thm}



There is a long history about the consequences of Bohr's inequality, in particular.
Indeed, Bohr \cite{Bohr-14} discovered that $\textbf{B}\geq 1/6$ and
the fact that $\textbf{B}=1/3$ was obtained independently by M. Riesz, I. Schur and N. Weiner.
Extensions and modifications of Bohr's result can be found from
\cite{bib:22,PPS2002,bib:09} and the recent articles
\cite{bib:18,bib:30,BhowDas-18,
bib:21,bib:24,bib:25,bib:26,bib:31}). We refer to
\cite{bib:11,bib:12,bib:13,bib:14} for the extension of the Bohr
inequality to several complex variables. More recently, Kayumov and
Ponnusamy \cite{bib:16} introduced and investigated
Bohr--Rogosinski's radii for the family $\mathcal B$, and they
discussed Bohr--Rogosinski's radius for the class of subordinations.
In \cite{AiElSh-05}, Aizenberg, et al. generalized the Rogosinski
radius for holomorphic mappings of the open unit polydisk into an
arbitrary convex domain. In \cite{bib:15},  Kayumov et al.
investigated Bohr's radius for complex-valued harmonic mappings that
are locally univalent in $\ID$. Several improved versions of Bohr's
inequality were given by Kayumov and Ponnusamy in \cite{bib:16} and
these were subsequently followed by  Evdoridis et al.  \cite{bib:19}
to obtain improved versions of Bohr's inequality for the class of
harmonic mappings. In \cite{bib:17}, Kayumov and Ponnusamy discussed
Bohr's radius for the class of analytic functions $g$, when $g$ is
subordinate to a member of the class of odd univalent functions. For
more information about Bohr's inequality and further  related works,
we refer the reader to the recent survey article \cite{AAP2016} and
the references therein.

In this paper we shall introduce and investigate another version of
the Rogosinski inequality for  analytic functions defined on the
unit disk $\mathbb{D}$ by substituting the derivatives of the
analytic function instead of the coefficients of its power series.
We shall also introduce and study several new versions of the
classical Bohr's inequality.

\section{An improved version of the classical Rogosinski inequality}
What could happened to the partial sums of the analytic function in
the unit  disk if we replaced the coefficients
$a_0,a_1,\ldots,a_{N-1}$ by the functions
$f(z),f'(z),\ldots,f^{(N-1)}(z)$? In this section we give an answer
in the following form.

\bthm\label{AKS2-th1}
  Suppose that $f\in {\mathcal B}$ and $f(z)= \sum_{n=0}^\infty a_n z^n$. Then
\begin{equation*}
\left|\sum_{k=0}^n\frac{f^{(k)}(z)}{k!}z^k\right|\le \sum_{k=0}^n \binom{-\frac12}{k}^2 ~\mbox{ for all $\ds |z|\leq r\le\frac12$}.
\end{equation*}
\ethm
\bpf
To prove this theorem we will use a modification of Landau's method (see \cite{bib:32} and \cite[p.~26]{LanDai-86}).

We consider the function $g:\, \ID \to \ID$ defined by $g(\zeta)=f(\alpha (\zeta +1))$, where $|\alpha|\leq 1/2$, and use the
substitution $\xi= D(\zeta) = \alpha (\zeta +1)$. In view of the Cauchy integral formula, integration
along a circle $\gamma$ around the origin lying in its neighborhood, we have
\begin{equation*}
\frac{f^{(k)}(\alpha)}{k!}  = \frac1{2\pi i}\int_{D(\gamma)} \frac{f(\xi)}{(\xi -\alpha)^{k+1}}\,d\xi
= \frac1{2\pi i \alpha^k}\int_{\gamma } \frac{g(\zeta)}{\zeta^{k+1}}\,d\zeta
\end{equation*}
and thus,  we can write
\begin{equation}\label{inqty_3}
\sum_{k=0}^n\frac{f^{(k)}(\alpha)}{k!}\alpha^k = \frac1{2\pi i}\int_{\gamma} g(\zeta) \left (\sum_{k=0}^n \frac{1}{\zeta^{k+1}}\right )
 d\zeta = \frac1{2\pi i}\int_{\gamma} \frac{g(\zeta)}{\zeta^{n+1}} \left (\sum_{k=0}^n \zeta^k\right ) d\zeta.
\end{equation}
Set
$1+\zeta+\zeta^2+\zeta^3+\dots = (1-\zeta)^{-1}  =  K^2(\zeta)= (K_n(\zeta))^2 +O(\zeta ^{n+1}),
$
where we write
$$
K(\zeta) =  (1-\zeta)^{-1/2}  =\sum_{k=0}^\infty \binom{-\frac12}{k} (-\zeta)^k ~\mbox{ and }~
K_n(\zeta) = \sum_{k=0}^n \binom{-\frac12}{k} (-\zeta)^k.
$$

In view of the above observations, \eqref{inqty_3} reduces to
\begin{equation}\label{eqty_1}
\sum_{k=0}^n\frac{f^{(k)}(\alpha)}{k!}\alpha^k = \frac1{2\pi i}\int_{\gamma} \frac{g(\zeta)}{\zeta^{n+1}} (K_n(\zeta))^2 \,d\zeta
\end{equation}
and therefore, with $\zeta=|\zeta|e^{i\phi}$, and $|g(\zeta)|\le1$ for all $|\alpha|\le1/2$ and $|\zeta|\le1$, we have
\begin{eqnarray*}
\left|\sum_{k=0}^n\frac{f^{(k)}(\alpha)}{k!}\alpha^k\right|
&\le& \frac1{2\pi}\int_0^{2\pi} \frac{1}{|\zeta|^{n+1}} |K_n(\zeta)|^2|\zeta| \,d\phi
\, =\,  \frac1{|\zeta|^n} \sum_{k=0}^n \binom{-\frac12}{k}^2 |\zeta|^{2k}.
\end{eqnarray*}
Allowing $|\zeta|\ra 1$, we get
\begin{equation*}
\left|\sum_{k=0}^n\frac{f^{(k)}(\alpha)}{k!}\alpha^k\right|
\le   \sum_{k=0}^n \binom{-\frac12}{k}^2 ~\mbox{ for all $ \ds |\alpha|\le \frac12$}
\end{equation*}
which completes the proof of the theorem.
\epf

\section{Improved versions of the classical Bohr's inequality}

For $f\in {\mathcal B}$ and $f(z)= \sum_{n=0}^\infty a_n z^n$, the following inequalities due to Schwarz-Pick will be used frequently: for $|z|=r<1$,
\begin{equation}\label{AKP-SchPick}
|f(z)|\le \frac{r+a}{1+ar} ~\mbox{ and }~ |f'(z)|\le \frac{1-|f(z)|^2}{1-|z|^2},
\end{equation}
where $|a_0|=a\in [0,1)$. Also, it is well-known that the Taylor coefficients of  $f\in {\mathcal B}$ satisfy the inequalities:
\begin{equation}\label{AKP-SchPick2}
|a_k|\le 1-a^2 ~\mbox{ for each $k\ge 1$.}
\end{equation}
More generally, we have (\cite{bib:29}) the sharp estimate
\be\label{AKPs-Rusheq2}
\frac{|f^{(k)}(z)|}{k!}\le \frac{1-|f(z)|^2}{(1-|z|)^k(1+|z|)} ~\mbox{ for each $k\ge 1$ and $z\in\ID$,}
\ee
which in particular gives second inequality in \eqref{AKP-SchPick}, and \eqref{AKP-SchPick2} by setting
$z=0$ in \eqref{AKPs-Rusheq2}. In the following we also assume that $m\in \IN$, and the idea of
replacing $a_k$ by $\frac{f^{(k)}(z)}{k!}$ is used in \cite{bib:31}. But our concern here is slightly different
from theirs.

\bthm\label{AKP2-th1}
Suppose that $f\in {\mathcal B}$ and $f(z)= \sum_{n=0}^\infty a_n z^n$. Then
\begin{equation}\label{AKP-eq1}
A_f(z):=|f(z^m)|+|z^m|\,|f'(z^m)| +\sum_{k=2}^\infty |a_k|r^k \le 1~\mbox{ for all }~ r \le R_{m,1},
\end{equation}
where $R_{m,1}$ is the maximal positive root of the equation $\varphi_m(r)=0$ with
\begin{equation}\label{AKP-eq2}
\varphi_m(r)=(1-r)(r^{2m}+2r^m-1)+2r^2(1+r^m)^2
\end{equation}
and the constant $R_{m,1}$ cannot be improved.
\ethm
\begin{table}[!htbp]
\begin{center}
\begin{tabular}{|l|l|}
\hline
 {\bf $m$}   & $R_{m,1}$ \\
  \hline
1  & 0.280776   \\
  \hline
2  & 0.39149   \\
  \hline
3  & 0.441112   \\
  \hline
4  & 0.467644   \\
  \hline
5  & 0.482442   \\
  \hline
\end{tabular}
\end{center}
\caption{$R_{m,1}$ is the maximal positive root of the equation
$(1-r) (r^{2m}+2r^m-1) + 2r^2(1+r^m)^2=0$\label{AKP2-th1-T1}}
\end{table}
\bpf
Let $f\in {\mathcal B}$ and $|a_0|=a\in [0,1)$. 
 It is a simple exercise to see that for $0\leq x\leq x_0 ~(\leq 1)$ and $0\leq \alpha \leq 1/2$, we have
$b(x):=x+\alpha (1-x^2)\leq b(x_0). $ This simple fact will be used
in the later theorems also. Using this inequality and
\eqref{AKP-SchPick2}, we easily obtain from \eqref{AKP-SchPick} and
\eqref{AKP-eq1}  that
\begin{align*}
A_f(z) &\le |f(z^m)| + \frac{r^m}{1-r^{2m}}(1-|f(z^m)|^2)
+(1-a^2)\frac{r^2}{1-r}\\
& \le \frac{r^m+a}{1+ar^m} +
\frac{r^m}{1-r^{2m}}\left[1-\left(\frac{r^m+a}{1+ar^m}\right)^2\right]
+(1-a^2)\frac{r^2}{1-r}\\
&= \frac{r^m+a}{1+ar^m} + (1-a^2)\frac{r^m}{(1+ar^m)^2} +(1-a^2)\frac{r^2}{1-r}\\
&= 1-\frac{(1-a)(1-r^m)}{1+ar^m} + (1-a^2)\frac{r^m}{(1+ar^m)^2}
+(1-a^2)\frac{r^2}{1-r}
=1+\frac{(1-a)\varphi_m(a,r)}{(1+ar^m)^2(1-r)},
\end{align*}
where
\begin{align*}
\varphi_m(a,r) &= -(1-r^m)(1+ar^m)(1-r) +(1+a)r^m(1-r) +
r^2(1+a)(1+ar^m)^2 \\
&= (1-r)(ar^{2m}+2r^m-1) + r^2(1+a)(1+ar^m)^2.
\end{align*}
The second inequality above is justified because of the fact that
$\frac{r^m}{1-r^{2m}}\leq \frac{1}{2}$ for $r \leq
\sqrt[m]{\sqrt2-1}.$ Also,  $R_{m,1}\leq \sqrt[m]{\sqrt2-1}$, where
$R_{m,1}$ is as in the statement. Now, since $\varphi_m(a,r) $ is an
increasing function of $a$ in $[0,1)$, it follows that
$$\varphi_m(a,r)\leq \varphi_m(1,r)=(1-r)(r^{2m}+2r^m-1) + 2r^2(1+r^m)^2= \varphi_m(r),
$$
where $\varphi_m(r)$ is given by (\ref{AKP-eq2}). Clearly,
$A_f(z)\leq 1$ if $\varphi_m(r)\le 0$, which holds for $r\le
R_{m,1}$.

To show the sharpness of the radius $R_{m,1}$, we let $a \in[0,1)$ and consider the function
\begin{equation}\label{AKP2-eqex1}
f(z)=\frac{a+z}{1+az}=a+(1-a^2)\sum_{k=1}^\infty (-a)^{k-1} z^k, \ \ z\in \ID
\end{equation}
so that
$$\frac{f^{(k)}(z)}{k!}=(1-a^2)\frac{ (-a)^{k-1}}{(1+az)^{k+1}} ~\mbox{ for  $k\ge 1$ and $z\in\ID$.}
$$
For this function, we observe that for $z=r$ and $a\in [0,1)$,
\beq\label{AKP-eq3}
|f(z^m)|+|z^m|\,|f'(z^m)| +\sum_{k=2}^\infty |a_k|r^k
&=& \frac{a+2r^m+ar^{2m}}{(1+ar^m)^2}+(1-a^2)\frac{ar^2}{1-ar}\nonumber\\
&=& 1+\frac{(1-a)P_m(a,r) }{(1+ar^m)^2(1-ar)}, \eeq where
$P_m(a,r)=(1-ar)(ar^{2m}+2r^m-1) + ar^2(1+a)(1+ar^m)^2
$
and the last expression \eqref{AKP-eq3} is larger than $1$ if and
only if $P_m(a,r)>0$. By a simple calculation, we find that
$$\frac{\partial P_m(a,r)}{\partial a}= r^{2m} +2a(r^2-r^{2m+1}) +r(1+r-2r^m)+3a^2r^{2m+2}+4ar^{m+2}+2ar^2+4a^3r^{2m+2}+6a^2r^{m+2}
$$
which is clearly non-negative for each $r\in [0,1)$ and thus, for each $r\in [0,1)$, $P_m(a,r)$ is an increasing function of $a$. This
fact gives
$$0< r^{2m} +r(1+r-2r^m) =P_m(0,r)\leq P_m(a,r)\le P_m(1,r)= \varphi_m(r),
$$
where $\varphi_m(r)$ is given by (\ref{AKP-eq2}).
Therefore, the right hand side of (\ref{AKP-eq3}) is smaller than or equal to $1$ for all $a\in [0,1)$, only in the case $r\le R_{m,1}$.
Finally, it also suggests that the right hand side of (\ref{AKP-eq3}) is larger than $1$ if $r > R_{m,1}$. This completes the proof.
\epf

\br
In Tables \ref{AKP2-th1-T1} and \ref{AKP2-th2-T2}, we listed the values of $R_{m,1}$ and $R_{m,2}$ for certain values of $m$.
If we allow $m\to\infty$ in Theorem \ref{AKP2-th1} (resp. Theorem \ref{AKP2-th2} below),
we note that $R_{m,1}\to 1/2$ (resp. $R_{m,2}\to 1/2$ below) and thus if $f\in {\mathcal B}$, then we have the
inequality
$$
|f(0)|+\sum_{k=2}^\infty \left|\frac{f^{(k)}(0)}{k!}\right|r^k \le 1 ~\mbox{ for all $ r\le 1/2$,}
$$
and the number $1/2$ is sharp.
\er

\bthm\label{AKP2-th2}
If $f\in {\mathcal B}$, then
\begin{equation}\label{AKP-eq5}
B_f(z):= |f(z^m)| + \sum_{k=2}^\infty \left|\frac{f^{(k)}(z^m)}{k!}\right|r^k \le 1~\mbox{for all $r \le R_{m,2}$},
\end{equation}
where $R_{m,2}$ is the minimum positive root of the equation $\psi _m(r)=0$ with
\begin{equation}\label{AKP-eq6}
\psi _m(r)=2r^2-(1-r^{2m})(1-r^m-r)
\end{equation}
and the constant $R_{m,2}$ cannot be improved.
\ethm
\begin{table}[!htbp]
\begin{center}
\begin{tabular}{|l|l|}
\hline
 {\bf $m$}   & $R_{m,2}$ \\
  \hline
1  & 0.355416   \\
  \hline
2  & 0.430586   \\
  \hline
3  & 0.464327 \\
  \hline
4  &  0.481418 \\
  \hline
5  & 0.490359  \\
  \hline
\end{tabular}
\end{center}
\caption{$R_{m,2}$ is the maximal positive root of the equation
$2r^2-(1-r^{2m})(1-r^m-r)=0$\label{AKP2-th2-T2}}
\end{table}
\bpf
As before we let $f\in {\mathcal B}$ and $a=|a_0|$.
By assumption, \eqref{AKP-SchPick} and \eqref{AKPs-Rusheq2} (with $z^m$ in place of $z$), we have
\begin{eqnarray*}
B_f(z)&\le& |f(z^m)|+ \frac{1-|f(z^m)|^2}{1+r^m}\sum_{k=2}^\infty \left (\frac{r}{1-r^m}\right )^k \\
&= & |f(z^m)|+ \frac{r^2}{(1-r^{2m})(1-r^m-r)}(1-|f(z^m)|^2)\\
&\le& \frac{r^m+a}{1+ar^m} + \frac{r^2}{(1-r^{2m})(1-r^m-r)}\left[1-\left(\frac{r^m+a}{1+ar^m}\right)^2\right] \\
&=& 1-\frac{(1-a)(1-r^m)}{1+ar^m} +
\frac{r^2(1-a^2)}{(1+ar^m)^2(1-r^m-r)}=
1+\frac{(1-a)\psi_m(a,r)}{(1+ar^m)^2(1-r^m-r)},
\end{eqnarray*}
where
\beqq
\psi_m(a,r) &=& -(1-r^m)(1+ar^m)(1-r^m-r)+r^2(1+a).
\eeqq
The second inequality above is a consequence of our earlier observation used in Theorem \ref{AKP2-th1} but this time with
$\alpha = r^2/[(1-r^{2m})(1-r^m-r)]$. It is a simple exercise to see that $\psi_m(a,r)$, for each $m\geq 1$,
is an increasing function of $a$ in $[0,1)$, and thus, it follows that
$
\psi_m(a,r) \leq  \psi_m(1,r) =\psi_m(r),
$
where $\psi_m(r)$ is given by \eqref{AKP-eq6}.
Clearly, $\psi_m(a,r)\le0$ if $\psi_m(r)\le 0$, which holds for $r\le R_{m,2}$, where
$R_{m,2}$ is the minimum positive root of the equation $\psi_m(r)= 0$. Thus, $B_f(z)\leq 1$
for $r\le R_{m,2}$ and the inequality \eqref{AKP-eq5} follows.

To show the sharpness of the radius $R_{m,2}$, we let $a \in(0,1)$ and consider the function
$g(z)=f(-z)$, where $f$ is given by \eqref{AKP2-eqex1}. Then for $g$, we easily have
$$\frac{g^{(k)}(z)}{k!} =-(1-|a|^2)\frac{a^{k-1}}{(1-az)^{k+1}},  \ \ z\in \ID.
$$
Now, we choose $a$ as close to $1$ as we please and set $z=r<\sqrt[m]a$. By a simple calculation,
the corresponding $B_g(z)$  takes the form
\beq\label{AKP-eq7}
\hspace{-.8cm}B_g(z) &=& \frac{a-r^m}{1-ar^m} + \frac{ar^2(1-a^2)}{(1-ar^m)^2(1-ar^m-ar)}\nonumber\\
&=& 1-\frac{(1-a)(1+r^m)}{1-ar^m} +
\frac{ar^2(1-a^2)}{(1-ar^m)^2(1-ar^m-ar)}=
1+\frac{(1-a)P_m(a,r)}{(1-ar^m)^2(1-ar^m-ar)},
\eeq
where
$$P_m(a,r)=a(1+a)r^2-(1+r^m)(1-ar^m)(1-ar^m-ar).
$$
Clearly, $B_g(z)<1$ if and only if $P_m(a,r)<0$. By elementary calculations, we find that
\beqq
\frac{\partial P_m(a,r)}{\partial a}&=& 2a[r^2-r^m(1+r^{m})(r^{m}+r)]+r^2+(1+r^{m})(2r^m+r)\\
&=&2a(r^2-r^{2m}-r^{m+1}-r^{3m}-r^{2m+1})+(2r^{2m}+r^{m+1}+2r^m+r^2+r)
\eeqq
which is easily seen to be greater than or equal to $0$ for any $r\in [0,1)$ and $m\ge 1$. Consequently,
$$P_m(a,r)\le P_m(1,r)=\psi _m(r) =2r^2-(1-r^{2m})(1-r^m-r).
$$
Therefore, the expression on the right of (\ref{AKP-eq7}) is smaller than or equal to $1$ for all $a\in (0,1)$, only in the case when $r\le R_{m,2}$.
Finally, it also suggests that $a\to 1$ in the right hand side of (\ref{AKP-eq7}) shows that  the expression (\ref{AKP-eq7}) is larger than
$1$ if and only if $r > R_{m,2}$. This completes the proof.
\epf

\bthm\label{AKP2-th3}
Suppose that $f\in {\mathcal B}$ and $f(z)= \sum_{n=0}^\infty a_n z^n$. Then
\begin{equation}\label{AKP-eq15}
C_f(z):=|f(z^m)|+|z|\,|f'(z^m)| +\sum_{k=2}^\infty |a_k|r^k \le 1\quad\text{for all}\quad r \le R_{m,3},
\end{equation}
where $R_{m,3}$ is the maximal positive root of the equation $\Phi_m (r)=0$ with
\begin{equation}\label{AKP-eq16}
\Phi_m (r)= 3r-1 + r^m\left[2r^2(r^m+2)+r^m(1-r)\right]
\end{equation}
and the constant $R_{m,3}$ cannot be improved.
\ethm
\begin{table}[!htbp]
\begin{center}
\begin{tabular}{|l|l|}
\hline
 {\bf $m$}   & $R_{m,3}$ \\
  \hline
1  & 0.280776   \\
  \hline
2  & 0.316912   \\
  \hline
3  & 0.327911 \\
  \hline
4  & 0.33152 \\
  \hline
5  & 0.332726  \\
  \hline
\end{tabular}
\end{center}
\caption{$R_{m,3}$ is the maximal positive root of the equation
$3r-1 + r^m[2r^2(r^m+2)+r^m(1-r)]=0$ \label{AKP2-th3-T3}}
\end{table}
\bpf
As in the proofs of Theorems \ref{AKP2-th1} and \ref{AKP2-th2}, it follows from \eqref{AKP-SchPick}, \eqref{AKP-SchPick2} and \eqref{AKPs-Rusheq2} that
\begin{eqnarray*}
C_f(z) &\le& |f(z^m)| + \frac{r}{1-r^{2m}}(1-|f(z^m)|^2) +(1-a^2)\frac{r^2}{1-r} \\
&\le& \frac{r^m+a}{1+ar^m} + \frac{r}{1-r^{2m}}\left[1-\left(\frac{r^m+a}{1+ar^m}\right)^2\right] +(1-a^2)\frac{r^2}{1-r} \\
&=& 1-\frac{(1-a)(1-r^m)}{1+ar^m} + (1-a^2)\frac{r}{(1+ar^m)^2} +(1-a^2)\frac{r^2}{1-r} \\
&=& \frac{a+2r^m+ar^{2m}}{(1+ar^m)^2}+(1-a^2)\frac{r^2}{1-r}=1+\frac{(1-a)\Phi_m(a,r)}{(1+ar^m)^2(1-r)}
\end{eqnarray*}
where
$$
\Phi_m(a,r) =
-(1-r^m)(1+ar^m)(1-r) + r(1+a)(1-r) +  r^2(1+a)(1+ar^m)^2
$$
$$
= r(1+a) + ar^{m+2}(1+a)(2+ar^m)-(1-r^m)(1+ar^m)(1-r)
$$
$$
\le \Phi_m(1,r) =\Phi_m(r) 
$$
because $\Phi_m(a,r) $ is seen to be an increasing function of $a$ in
$[0,1)$, and $\Phi_m(r)$ is given by \eqref{AKP-eq16}. Note that the
second inequality above holds since $\max \frac{2r}{1-r^{2m}}<1$ and
so for any $r <R_m$, where $R_m$ is the maximal positive root of the
equation $2r-(1-r^{2m})=0$, and $R_{m,3}<R_m$ for $m\in \IN$, where
$R_{m,3}$ is the maximal positive root of the equation $\Phi_m
(r)=0$. Since $\Phi_m(r)\le0$  for $r\le R_{m,3}$, we obtain
$C_f(z)\leq 1$ for $r\le R_{m,3}$ and the desired inequality
\eqref{AKP-eq15} follows.

It remains to show the sharpness of the radius $R_{m,3}$. To do this we let $a \in[0,1)$ and consider the function
$f$ is given by \eqref{AKP2-eqex1}.
For this function, we observe that for $z=r$,
\begin{equation}\label{AKP-eq17}
C_f(z) = \frac{(r^m+a)(1+ar^m)+r(1-a^2)}{(1+ar^m)^2}+(1-a^2)\frac{ar^2}{1-ar}= 1+\frac{(1-a)Q_m(a,r)}{(1+ar^m)^2(1-r)},
\end{equation}
where
$$Q_m(a,r)=r(1+a) + a^2r^{m+2}(1+a)(2+ar^m)-(1-r^m)(1+ar^m)(1-ar).
$$
We see that   $C_f(z)>1$ for $a \in[0,1)$ if and only if $Q_m(a,r) >0$.
By a simple calculation, we find that $Q_m(a,r)$ is an increasing function of $a$ in $[0,1)$
and therefore, we have
\begin{eqnarray*}
Q_m(a,r)\le Q_m(1,r) = 2r + 2r^{m+2}(2+r^m)-(1-r^m)(1+r^m)(1-r)=\Phi_m(r),
\end{eqnarray*}
where $\Phi_m(r)$ is given by \eqref{AKP-eq16}.
Therefore, the expression (\ref{AKP-eq17}) is smaller than or equal to $1$ for all $a\in [0,1)$, only  when
$r\le R_{m,3}$. Finally, it also suggests that $a\to 1$ in (\ref{AKP-eq17}) shows that the expression (\ref{AKP-eq17}) is
larger than $1$ if $r > R_{m,3}$. This completes the proof.
\epf

\br
In Table \ref{AKP2-th3-T3}, we listed the values of $R_{m,3}$ for certain values of $m$.
If we allow $m\to\infty$ in Theorem \ref{AKP2-th3}, we see that $R_{m,3}\to \frac13$ and hence we have the classical Bohr inequality
for $f\in {\mathcal B}$:
$$ |f(0)| +\sum_{k=1}^\infty |a_k|\,|z|^k \le 1 ~\mbox{ for all $r\le1/3$},
$$
and $1/3$ is sharp.
\er


\section{Two Improved versions of Bohr's inequality for harmonic mappings}

\bthm\label{AKP2-th4}
Suppose that $f(z)= h(z) + \overline{g(z)} =\sum_{n = 0}^\infty a_nz^n
+ \overline{\sum_{n = 1}^\infty b_nz^n}$ is a harmonic mapping of $\ID$ such that
$|g'(z)|\leq k|h'(z)|$ for some $k\in [0,1]$ and $h\in {\mathcal B}$. Then we have
\begin{equation}\label{AKP-eq9}
D_f(z):=
|h(z^m)|+\sum_{n=1}^\infty |a_n|r^n + \sum_{n=1}^\infty |b_n|r^n \le 1~\mbox{ for all $r \le R^k_{m,1}$},
\end{equation}
where $R^k_{m,1}$ is the maximal positive root of the equation $\lambda_m (r)=0$ with
\begin{equation}\label{AKP-eq10}
\lambda_m (r)= 2r(1+k)(1+r^m)-(1-r)(1-r^m)
\end{equation}
and the constant $R^k_{m,1}$ cannot be improved.
\ethm
\begin{table}[!htbp]
\begin{center}
\begin{tabular}{|l|l|}
\hline
 {\bf $m$}   & $R^1_{m,1}$ \\
  \hline
1  &  0.154701  \\
  \hline
2  & 0.188829   \\
  \hline
3  & 0.197544 \\
  \hline
4  & 0.199494 \\
  \hline
5  & 0.199898  \\
  \hline
\end{tabular}
\end{center}
\caption{$R^1_{m,1}$ is the maximal positive root of the equation
$4r(1+r^m)-(1-r)(1-r^m)=0$ \label{AKP2-th4-T4}}
\end{table}
\bpf
Recall that, as $h\in {\mathcal B}$ and $h(z)=\sum_{n = 0}^\infty a_nz^n$,
$$h(z)\prec  \frac{a_0+z}{1+\overline{a_0}z}=a_0+(1-|a_0|^2)\sum_{n=1}^\infty (-1)^{n-1}(\overline{a_0})^{n-1}z^n, \ \ z\in \ID,
$$
which gives \cite{BhowDas-18,bib:30}
\begin{eqnarray*}
\sum_{n=1}^{\infty}|a_{n}|r^{n}\leq (1-|a_0|^2)\sum_{n=1}^{\infty} |a_0|^{n-1}r^{n}= \frac{(1-a^2)r}{1-ar}~\mbox{ for all }\, r \leq\frac{1}{3},
\end{eqnarray*}
where $a=|a_0|\in[0,1)$. Moreover, by assumption, we obtain that $g'(z)\prec_q k h'(z)$ which quickly gives from  \cite{bib:30} that
\begin{eqnarray*}
\sum_{n=1}^{\infty}n|b_{n}|r^{n-1}\leq \sum_{n=1}^{\infty}k n\, |a_{n}|r^{n-1}~\mbox{ for all }\, r \leq\frac{1}{3}
\end{eqnarray*}
and integrating this with respect to $r$ gives
$$\sum_{n=1}^{\infty}|b_{n}|r^{n}\leq k\sum_{n=1}^{\infty}|a_{n}|r^{n}~\mbox{ for all }\, r \leq\frac{1}{3}.
$$
Here $\prec _q$ denotes the quasi-subordination. Using these and the first inequality in \eqref{AKP-SchPick} for $h(z)$ one can
obtain that for $|z|=r\le 1/3$,
\begin{align*}
D_f(z) &\le  \frac{r^m+a}{1+ar^m} + (1+k)r\frac{1-a^2}{1-ar}
=1+\frac{(1-a)\lambda_m(a,r)}{(1+ar^m)(1-ar)},
\end{align*}
where
$$\lambda_m(a,r) = r(1+k)(1+a)(1+ar^m)-(1-r^m)(1-ar),
$$
which is indeed an  increasing function of $a\in [0,1)$ so that
$\lambda_m(a,r) \leq \lambda_m(1,r)= \lambda_m(r),$ where $\lambda_m (r)$
is given by \eqref{AKP-eq10}. We see that $D_f(z)\leq 1$ if $\lambda
_m(r)\le0$, which holds for $r\le R^k_{m,1}$, where $R^k_{m,1}$ is
the maximal positive root of the equation  $\lambda_m (r)=0$. This
proves the inequality \eqref{AKP-eq9}.

Finally, to show the sharpness of the radius $R^k_{m,1}$,  we consider the function
\be\label{AKP-eq18}
f(z)=h(z) +\lambda k \overline{h(z)}, \quad h(z)= \frac{z+a}{1+az},
\ee
where  $\lambda \in (0,1]$. For this function, we get that (for $z=r$ and $a\in [0,1)$)
$$ D_f(z)  = \frac{r^m+a}{1+ar^m} + (\lambda k +1) r\frac{1-a^2}{1-ra}
$$
and the last expression shows the sharpness of $R^k_{m,1}$ with $\lambda \to 1$. This completes the proof of the theorem.
\epf

\br In Table \ref{AKP2-th4-T4}, we listed the values of $R^k_{m,1}$ for $k=1$ and for certain values of $m$.
When $m\to\infty$, we have from Theorem \ref{AKP2-th4} that $R^k_{m,1}\to \frac{2}{4k+6}$. Thus, under the hypotheses
of Theorem \ref{AKP2-th4}, we have
$$|h(0)|+\sum_{n=1}^\infty |a_n|r^n + \sum_{n=1}^\infty |b_n|r^n \le 1~\mbox{ for all $\ds r \le \frac{2}{4k+6}$}
$$
which for $k=0$ gives the classical Bohr's inequality and for $k=1$, this inequality contains the Bohr inequality for
sense-preserving harmonic mapping  $f(z)= h(z) + \overline{g(z)}$ of the disk $\ID$ with the Bohr radius $1/5$
(see \cite{bib:15}).
\er

\bthm\label{AKP2-th5}
Assume the hypotheses of Theorem \ref{AKP2-th4}. Then we have
$$
E_f(z):=
|h(z^m)|+|z^m|\,|h'(z^m)|+\sum_{n=2}^\infty |a_n|r^n + \sum_{n=1}^\infty |b_n|r^n \le 1 ~\mbox{ for all $r \le R^k_{m,2}$},
$$
where $R^k_{m,2}$ is the maximal positive root of the equation $\Lambda_m (r)=0$ with
$$
\Lambda_m (r)=
(1-r)(r^{2m}+2r^m-1)+2r(r+k)(1+r^m)^2
$$
and the constant $R^k_{m,2}$ cannot be improved.
\ethm
\begin{table}[!htbp] \begin{center}
\begin{tabular}{|l|l|}
\hline
 {\bf $m$}   & $R^1_{m,2}$ \\
  \hline
1  &  0.1671  \\
  \hline
2  &  0.240751  \\
  \hline
3  & 0.267472 \\
  \hline
4  & 0.276691\\
  \hline
5  & 0.279585  \\
  \hline
\end{tabular}
\end{center}
\caption{$R^1_{m,2}$ is the maximal positive root of the equation
$(1-r)(r^{2m}+2r^m-1)+2r(r+1)(1+r^m)^2=0$ \label{AKP2-th5-T5}}
\end{table}
\bpf
As in the proofs of Theorem \ref{AKP2-th4} and earlier theorems, we easily have
\begin{align*}
E_f(z) &\le \frac{r^m+a}{1+ar^m} +
\frac{r^m}{1-r^{2m}}\left[1-\left(\frac{r^m+a}{1+ar^m}\right)^2\right]+
(1-a^2)\frac{r^2}{1-r} + k(1-a^2)\frac{r}{1-r}\\
& =  1-\frac{(1-a)(1-r^m)}{1+ar^m} +  \frac{(1-a^2)r^m}{(1+ar^m)^2} +\frac{(1-a^2)(r+k)r}{1-r}\\
&= 1+\frac{(1-a)\Lambda_m(a,r)}{(1+ar^m)^2(1-r)},
\end{align*}
where
\begin{eqnarray*}
\Lambda_m(a,r) &=&  - (1-r^m)(1-r)(1+ar^m) +  r^{m}(1+a)(1-r) + (1+a)r(r+k)(1+ar^m)^2 \\
&=& (1-r)(ar^{2m}+2r^m-1) + r(1+a)(r+k)(1+ar^m)^2\le \Lambda_m(1,r) 
=\Lambda_m(r).
\end{eqnarray*}
The first inequality above is justified with the same reasoning as in the proofs of earlier theorems.

Now, we see that $E_f(z)\leq 1$ whenever $\Lambda_m(r) \le0$, which holds for $r\le R^k_{m,2}$,
where $R^k_{m,2}$ is the maximal positive root of the equation $\Lambda_m(r)=0$.

To show the sharpness of the radius $R^k_{m,2}$,  consider the function $f$ defined by \eqref{AKP-eq18}
with $\lambda \in (0,1]$. For this function, the corresponding expression for $E_f(z)$ with $z=r$ turned out to be
\begin{equation}\label{AKP-eq13}
E_f(z)= \frac{a+2r^m+ar^{2m}}{(1+ar^m)^2}+\frac{ar(1-a^2)(r+\lambda)}{1-ar}.
\end{equation}
The last expression is larger than $1$ if and only if $P^k_m(a,r)>0$, where
\begin{equation}\label{AKP-eq14}
P^k_m(a,r)=(1-ar)(ar^{2m}+2r^m-1) + ar(1+a)(r+\lambda)(1+ar^m)^2.
\end{equation}
By a simple calculation, we find that $P^k_m(a,r)$ is an increasing function of $a\in [0,1)$, and for each $r\in [0,1)$, so that
$$P^k_m(a,r)\le P^k_m(1,r) = (1-r)(r^{2m}+2r^m-1) + 2r(r+\lambda)(1+r^m)^2.
$$
Therefore, the expression (\ref{AKP-eq13}) is smaller than or equal to $1$ for all $a\in [0,1)$, only in the case
when $r\le R^{k}_{m,2}$ ($\lambda =k$). Finally, it also suggests that $a\to 1$ in (\ref{AKP-eq14}) shows that  the expression
(\ref{AKP-eq13}) is larger than $1$ if $r > R^k_{m,2}$. This completes the proof.
\epf

\br In Table \ref{AKP2-th5-T5}, we listed the values of $R^k_{m,2}$ for $k=1$ and for certain values of $m$.
If we allow $m\to\infty$ in Theorem \ref{AKP2-th5}, we obtain that
$$R^k_{m,2}\to
R^k_{2}:=\frac14\left(\sqrt{(2k+1)^2+8}-(2k+1)\right),
$$
where $R^k_{2}$ is the positive root of the equation $2x(x+k)+x-1=0$
and the conclusion of Theorem~\ref{AKP2-th5} takes the following form:
$$ |h(0)|+\sum_{n=2}^\infty |a_n|r^n + \sum_{n=1}^\infty |b_n|r^n \le 1 ~\text{ for all }~
r \le \frac14\left(\sqrt{(2k+1)^2+8}-(2k+1)\right).
$$
\er

\subsection*{Acknowledgements}
The work of S. Alkhaleefah and I. Kayumov is supported by the Russian Science Foundation under grant 18-11-00115.
The  work of the third author is supported by Mathematical Research Impact Centric Support of Department of 
Science \& Technology, India  (MTR/2017/000367).


\end{document}